\documentclass[reqno]{amsart}

\usepackage[latin1]{inputenc}
\usepackage{amssymb}
\usepackage{graphicx}
\usepackage{amscd}
\usepackage[hidelinks]{hyperref}
\usepackage{color}
\usepackage{float}
\usepackage{graphics,amsmath,amssymb}
\usepackage{amsthm}
\usepackage{amsfonts}
\usepackage{latexsym}
\usepackage{epsf}
\usepackage{enumerate}
\usepackage{xifthen}
\usepackage{mathrsfs}
\usepackage{dsfont}
\usepackage{makecell}
\usepackage[FIGTOPCAP]{subfigure}
\usepackage{amsmath}
\allowdisplaybreaks[4]
\usepackage{listings}
\usepackage{etoolbox}
\usepackage{fancyhdr}

 \newtheoremstyle{mytheorem}
 {3pt}
 {3pt}
 {\slshape}
 {}
 {\bfseries}
 {.}
 { }
 {}

\numberwithin{equation}{section}

\theoremstyle{theorem}
\newtheorem{theorem}{Theorem}[section]

\newtheorem{lemma}[theorem]{Lemma}

\theoremstyle{definition}

\newcommand{\Keywords}[1]{\ifthenelse{\isempty{#1}}{}{\smallskip \smallskip \noindent \textbf{Keywords}. #1}}
\newcommand{\MSC}[2][2010]{\ifthenelse{\isempty{#2}}{}{\smallskip \smallskip \noindent \textbf{#1MSC}. #2}}
\newcommand{\abstractnote}[1]{\ifthenelse{\isempty{#1}}{}{\smallskip \smallskip \noindent \textsuperscript{\dag}#1}}

\makeatletter
\def\specialsection{\@startsection{section}{1}%
  \z@{\linespacing\@plus\linespacing}{.5\linespacing}%
  {\normalfont}}
\def\section{\@startsection{section}{1}%
  \z@{.7\linespacing\@plus\linespacing}{.5\linespacing}%
  {\normalfont\scshape}}
\patchcmd{\@settitle}{\uppercasenonmath\@title}{\Large\boldmath}{}{}
\patchcmd{\@settitle}{\begin{center}}{\begin{flushleft}}{}{}
\patchcmd{\@settitle}{\end{center}}{\end{flushleft}}{}{}
\patchcmd{\@setauthors}{\MakeUppercase}{\normalsize}{}{}
\patchcmd{\@setauthors}{\centering}{\raggedright}{}{}
\patchcmd{\section}{\scshape}{\large\bfseries\boldmath}{}{}
\patchcmd{\subsection}{\bfseries}{\bfseries\boldmath}{}{}
\renewcommand{\@secnumfont}{\bfseries}
\patchcmd{\@startsection}{\@afterindenttrue}{\@afterindentfalse}{}{}
\patchcmd{\abstract}{\leftmargin3pc}{\leftmargin1pc}{}{}

\def\maketitle{\par
  \@topnum\z@ 
  \@setcopyright
  \thispagestyle{empty}
  \ifx\@empty\shortauthors \let\shortauthors\shorttitle
  \else \andify\shortauthors
  \fi
  \@maketitle@hook
  \begingroup
  \@maketitle
  \toks@\@xp{\shortauthors}\@temptokena\@xp{\shorttitle}%
  \toks4{\def\\{ \ignorespaces}}
  \edef\@tempa{%
    \@nx\markboth{\the\toks4
      \@nx\MakeUppercase{\the\toks@}}{\the\@temptokena}}%
  \@tempa
  \endgroup
  \c@footnote\z@
  \@cleartopmattertags
}
\makeatother

\title[Congruences for 1-shell TSPPs]{Congruences for 1-shell totally symmetric plane partitions}

\author[S. Chern]{Shane Chern}
\address{Department of Mathematics, Pennsylvania State University, University Park, PA 16802, USA
\newline\indent
School of Mathematical Sciences, Zhejiang University, Hangzhou, 310027, China}
\email{shanechern@psu.edu}

\begin{document}

{\footnotesize\noindent \textit{Integers} \textbf{17} (2017), Paper No. A21, 7 pp.}

\bigskip \bigskip

\maketitle

\begin{abstract}  
Let $f(n)$ denote the number of 1-shell totally symmetric plane partitions of weight $n$. Recently,  Hirschhorn and Sellers, Yao, and Xia established a number of congruences modulo 2 and 5, 4 and 8, and 25 for $f(n)$, respectively. In this note, we shall prove several new congruences modulo 125 and 11 by using some results of modular forms. For example, for all $n\ge 0$, we have
\begin{align*}
f(1250n+125)&\equiv 0 \pmod{125},\\
f(1250n+1125)&\equiv 0 \pmod{125},\\
f(2750n+825)&\equiv 0 \pmod{11},\\
f(2750n+1925)&\equiv 0 \pmod{11}.
\end{align*}


\MSC{Primary 11P83; Secondary 05A17.}

%
%
%
\end{abstract}

\section{Introduction}

A plane partition is a two-dimensional array of integers $\pi_{i,j}$ that are weakly decreasing in both indices and that add up to the given number $n$, namely, $\pi_{i,j}\ge \pi_{i+1,j}$, $\pi_{i,j}\ge \pi_{i,j+1}$, and $\sum \pi_{i,j}=n$. If a plane partition is invariant under any permutation of the three axes, we call it a totally symmetric plane partition (see, e.g., Andrews \textit{et al.} \cite{APS05} and Stembridge \cite{Ste95} for more details). In 2012, Blecher \cite{Ble12} studied a special class of totally symmetric plane partitions which he called 1-shell totally symmetric plane partitions. A 1-shell totally symmetric plane partition has a self-conjugate first row/column (as an ordinary partition) and all other entries are 1. For example,
$$\begin{array}{cccc}
4 & 4 & 2 & 2 \\
4 & 1 & 1 & 1 \\
2 & 1 \\
2 & 1 \end{array}$$
is a totally symmetric plane partition.

Let $f(n)$ denote the number of 1-shell totally symmetric plane partitions of weight $n$, namely, the parts of the totally symmetric plane partition sum to $n$. In \cite{Ble12}, Blecher found the generating function of $f(n)$,
$$\sum_{n\ge 0}f(n)q^n=1+\sum_{n\ge 1}q^{3n-2}\prod_{i=0}^{n-2}\left(1+q^{6i+3}\right).$$
Recently, Hirschhorn and Sellers \cite{HS14}, Yao \cite{Yao14}, and Xia \cite{Xia15} established a number of congruences for $f(n)$, respectively. For example, for all $n\ge 0$, Hirschhorn and Sellers proved that
\begin{equation}\label{eq:HS}
f(10n+5)\equiv 0 \pmod{5},
\end{equation}
while Xia proved that
\begin{equation}
f(250n+125)\equiv 0 \pmod{25}.
\end{equation}
Moreover, Yao showed that, for all $n\ge 0$,
\begin{equation}
f(8n+3)\equiv 0 \pmod{4}.
\end{equation}

In this note, we shall prove several new congruences modulo 125 and 11 for $f(n)$. Here our methods are based on some results of modular forms, which are quite different from the proofs of the previous congruences. In fact, Radu and Sellers gave a strategy in \cite{RS11} to prove these Ramanujan-like congruences, and their methods can be tracked back to \cite{Ra09}. Our results are stated as follows.

\begin{theorem}\label{th:01}
For all $n\ge 0$, we have
\begin{equation}\label{eq:th1.1}
f(1250n+125)\equiv 0 \pmod{125}
\end{equation}
and
\begin{equation}\label{eq:th1.2}
f(1250n+1125)\equiv 0 \pmod{125}.
\end{equation}
\end{theorem}

\begin{theorem}\label{th:02}
For all $n\ge 0$, we have
\begin{equation}\label{eq:th2.1}
f(2750n+825)\equiv 0 \pmod{11}
\end{equation}
and
\begin{equation}\label{eq:th2.2}
f(2750n+1925)\equiv 0 \pmod{11}.
\end{equation}
\end{theorem}

By \eqref{eq:HS} and Theorem \ref{th:02}, we immediately get

\begin{theorem}\label{th:03}
For all $n\ge 0$, we have
\begin{equation}\label{eq:th3.1}
f(2750n+825)\equiv 0 \pmod{55}
\end{equation}
and
\begin{equation}\label{eq:th3.2}
f(2750n+1925)\equiv 0 \pmod{55}.
\end{equation}
\end{theorem}

\section{Preliminaries}

We first introduce some notations of \cite{RS11}. Let $M$ be a positive integer. We denote by $R(M)$ the set of integer sequences $\{r:r=(r_{\delta_1},\ldots,r_{\delta_k})\}$ indexed by the positive divisors $1=\delta_1<\cdots<\delta_k=M$ of $M$. For a positive integer $m$, let $[s]_m$ be the set of all elements congruent to $s$ modulo $m$. We also write $\mathbb{Z}_m^*$ the set of all invertible elements in $\mathbb{Z}_m$, and $\mathbb{S}_m$ the set of all squares in $\mathbb{Z}_m^*$. For $t\in\{0,\ldots,m-1\}$, we define by $\overline{\odot}_r$ the map $\mathbb{S}_{24m}\times\{0,\ldots,m-1\}$$\to$$\{0,\ldots,m-1\}$ with
$$([s]_{24m},t)\mapsto [s]_{24m}\overline{\odot}_r t\equiv ts+\frac{s-1}{24}\sum_{\delta\mid M}\delta r_{\delta}\pmod{m}.$$
Furthermore, we put $P_{m,r}(t):=\{[s]_{24m}\overline{\odot}_r t\ |\ [s]_{24m}\in\mathbb{S}_{24m}\}$.

Let $\Gamma:=SL_2(\mathbb{Z})$ and $\Gamma_\infty:=\left\{\left.\begin{pmatrix}1 & h \\0 & 1 \end{pmatrix}\ \right|\ h\in\mathbb{Z}\right\}$. For a positive integer $N$, we define the congruence subgroup of level $N$ as
$$\Gamma_0(N):=\left\{\left.\begin{pmatrix}a & b\\c & d\end{pmatrix}\in\Gamma\ \right|\ c\equiv 0\pmod{N}\right\}.$$
We also know that
$$[\Gamma:\Gamma_0(N)]=N\prod_{p\mid N}(1+p^{-1}),$$
where the product runs through the distinct prime numbers dividing $N$.

Now denote by $\Delta^*$ the set of tuples $(m,M,N,t,r=(r_\delta))$ which satisfy conditions given in \cite[p. 2255]{RS11}\footnote{According to a private communication between the author and S. Radu, the last condition of $\Delta^*$ should read: ``for $(s,j)=\pi(M,(r_\delta))$, if $2\mid m$, we have ($4\mid \kappa N$ and $8\mid Ns$) or ($2\mid s$ and $8\mid N(1-j)$).''}. Let $\kappa=\kappa(m)=\gcd(m^2-1,24)$. For $\gamma=\begin{pmatrix}a & b \\c & d \end{pmatrix}$, $r\in R(M)$, and $r'\in R(N)$, we set
$$p_{m,r}(\gamma)=\min_{\lambda\in\{0,\ldots,m-1\}}\frac{1}{24}\sum_{\delta\mid M}r_\delta\frac{\gcd^2(\delta(a+\kappa\lambda c),mc)}{\delta m}$$
and
$$p_{r'}^*(\gamma)=\frac{1}{24}\sum_{\delta\mid N} \frac{r'_\delta\gcd^2(\delta,c)}{\delta}.$$

Finally, we write $(a;q)_\infty :=\prod_{n\ge 0}(1-aq^n)$, and let
$$f_r(q):=\prod_{\delta\mid M}(q^{\delta};q^{\delta})_{\infty}^{r_\delta}=\sum_{n\ge 0}c_r(n)q^n$$
for some $r\in R(M)$. The following lemma (see \cite[Lemma 4.5]{Ra09} or \cite[Lemma 2.4]{RS11}) is a key to our proof.

\begin{lemma}\label{le:01}
Let $u$ be a positive integer, $(m,M,N,t,r=(r_\delta))\in\Delta^*$, $r'=(r'_\delta)\in R(N)$, $n$ be the number of double cosets in $\Gamma_0(N)\backslash\Gamma/\Gamma_\infty$ and $\{\gamma_1,\ldots,\gamma_n\}$ $\subset\Gamma$ be a complete set of representatives of the double coset $\Gamma_0(N)\backslash\Gamma/\Gamma_\infty$. Assume that $p_{m,r}(\gamma_i)+p_{r'}^*(\gamma_i)\ge 0$ for all $i=1,\ldots,n$. Let $t_{\min} := \min_{t'\in P_{m,r}(t)}t'$ and
$$v:=\frac{1}{24}\left(\left(\sum_{\delta\mid M}r_\delta+\sum_{\delta\mid N}r'_\delta\right)[\Gamma:\Gamma_0(N)]-\sum_{\delta\mid N}\delta r'_\delta\right)-\frac{1}{24m}\sum_{\delta\mid M}\delta r_\delta-\frac{t_{\min}}{m}.$$
Then if
$$\sum_{n=0}^{\lfloor v \rfloor}c_r(mn+t')q^n\equiv 0 \pmod{u}$$
for all $t'\in P_{m,r}(t)$, then
$$\sum_{n\ge 0}c_r(mn+t')q^n\equiv 0 \pmod{u}$$
for all $t'\in P_{m,r}(t)$.
\end{lemma}

\section{Proofs of the theorems}

\subsection{The upper bound}

In the first part of our proofs, we will compute the upper bound of $\lfloor v \rfloor$ in Lemma \ref{le:01} for each theorem. Let $g(n)$ be given by
\begin{equation}\label{eq:gn}
\sum_{n\ge 0}g(n)q^n:=\frac{(q^2;q^2)_{\infty}^3}{(q;q)_{\infty}^2}.
\end{equation}
In \cite{HS14}, Hirschhorn and Sellers proved that
\begin{equation}\label{eq:fg}
f(6n+1)=g(n).
\end{equation}
Moreover, we write
\begin{equation}\label{eq:gpa}
\sum_{n\ge 0}g_{\alpha,p}(n)q^n:=\frac{(q;q)_{\infty}^{p^\alpha-2}(q^2;q^2)_{\infty}^3}{(q^p;q^p)_{\infty}^{p^{\alpha-1}}},
\end{equation}
where $\alpha$ is a positive integer and $p$ is prime. By \cite[Lemma 1.2]{RS11}, we obtain
\begin{equation}\label{eq:ggpa}
\sum_{n\ge 0}g_{\alpha,p}(n)q^n\equiv \sum_{n\ge 0}g(n)q^n \pmod{p^\alpha}.
\end{equation}

Note that \cite[Theorem 2.1]{HS14} tells that $f(n)=0$ if $n\equiv 0,2$ (mod $3$) for all $n\ge 1$. We therefore have $f(1250\cdot 3n+125)=f(1250\cdot(3n+2)+125)=0$. To prove \eqref{eq:th1.1}, it suffices to prove $f(3750n+1375)=f(1250\cdot(3n+1)+125)\equiv 0$ (mod $125$), which yields
\begin{equation}\label{eq:g125.1}
g_{3,5}(625n+229)\equiv 0 \pmod{125}.
\end{equation}
Similarly, to prove \eqref{eq:th1.2}, \eqref{eq:th2.1}, and \eqref{eq:th2.2}, we only need to prove
\begin{equation}\label{eq:g125.2}
g_{3,5}(625n+604)\equiv 0 \pmod{125},
\end{equation}
\begin{equation}\label{eq:g11.1}
g_{1,11}(1375n+1054)\equiv 0 \pmod{11},
\end{equation}
and
\begin{equation}\label{eq:g11.2}
g_{1,11}(1375n+779)\equiv 0 \pmod{11},
\end{equation}
respectively.

Let
$$r^{(\alpha,p)}:=(r_1,r_2,r_p,r_{2p})=(p^\alpha-2,3,-p^{\alpha-1},0)\in R(2p).$$
By the definition of $P_{m,r}(t)$, we have
$$P_{m,r^{(\alpha,p)}}(t)=\left\{t'\ |\ t'\equiv ts+(s-1)/6\ (\bmod\ m),0\le t'\le m-1,[s]_{24m}\in\mathbb{S}_{24m}\right\}.$$
One readily verifies $P_{625,r^{(3,5)}}(229)=\{229,604\}$. Next we set
$$(m,M,N,t,r=(r_1,r_2,r_5,r_{10}))=(625,10,10,229,(123,3,-25,0))\in\Delta^*$$
and
$$r'=(r'_1,r'_2,r'_5,r'_{10})=(13,0,0,0).$$
Moreover, by \cite[Lemma 2.6]{RS11}, $\{\gamma_\delta:\delta\mid N\}$ contains a complete set of representatives of the double coset $\Gamma_0(N)\backslash\Gamma/\Gamma_\infty$ where $\gamma_\delta=\begin{pmatrix}
1 & 0 \\
\delta & 1 
\end{pmatrix}$. One may see that all these constants satisfy the assumption of Lemma \ref{le:01}. We therefore obtain $\lfloor v \rfloor=84$.

To get the upper bound $\lfloor v \rfloor$ for Theorem \ref{th:02}, we have $P_{1375,r^{(1,11)}}(1054)=\{779,1054\}$. Similarly, we can compute other relevant constants of \eqref{eq:g11.1} and \eqref{eq:g11.2}, which are listed in Table \ref{ta:01}.

\begin{table}[ht]
\caption{Relevant constants of \eqref{eq:g11.1} and \eqref{eq:g11.2}}\label{ta:01}
\renewcommand\arraystretch{1.25}
\noindent\[
\begin{tabular}{l}
\hline
$P_{1375,r^{(1,11)}}(1054)=\{779,1054\}$\\
$(m,M,N,t,r=(r_1,r_2,r_{11},r_{22}))=(1375,22,110,1054,(9,3,-1,0))$\\
$r'=(r'_1,r'_2,r'_5,r'_{10},r'_{11},r'_{22},r'_{55},r'_{110})=(6,0,0,0,0,0,0,0)$\\
$\lfloor v \rfloor=152$\\
\hline
\end{tabular}
\]
\end{table}

\subsection{Simplifying the verification}

We should notice that as $n$ approaches the upper bound $\lfloor v \rfloor$ in both theorems, the verification will become difficult. Hence we provide a method that can simplify the calculation. First, we notice that
\begin{equation}\label{eq:3.2.1}
\sum_{n\ge 0}g(n)q^n=\frac{(q^2;q^2)_{\infty}^3}{(q;q)_{\infty}^2}=\frac{1}{(q^2;q^2)_{\infty}}\left(\frac{(q^2;q^2)_{\infty}^2}{(q;q)_{\infty}}\right)^2.
\end{equation}
From \cite[Chapter 16, Entry 22(ii)]{Ber91} we know that
$$\frac{(q^2;q^2)_{\infty}^2}{(q;q)_{\infty}}=\sum_{n\ge 0}q^{T_n},$$
where $T_n=n(n+1)/2$ is the triangular number. Now we write
$$\left(\frac{(q^2;q^2)_{\infty}^2}{(q;q)_{\infty}}\right)^2=\sum_{n\ge 0}a(n)q^n.$$
Then
$$a(n)=\sharp\left\{(n_1,n_2)\in(\mathbb{N}\cup \{0\})^2: n=T_{n_1}+T_{n_2}\right\}.$$
Notice that $n=T_{n_1}+T_{n_2}$ implies $8n+2=(2n_1+1)^2+(2n_2+1)^2$. Now if $8n+2$ is a sum of two squares, then both the squares are odd. This is because a square is congruent to $0,1,4$ modulo $8$. We therefore have
$$4a(n)=r(8n+2),$$
where $r(n)$ denotes the number of representations of $n$ by two squares. For example, $r(5)=8$ since $5=(\pm 1)^2+(\pm 2)^2=(\pm 2)^2+(\pm 1)^2$.

Let $p(n)$ be the partition function given by
$$\sum_{n\ge 0}p(n)q^n:=\frac{1}{(q;q)_{\infty}}.$$
It follows by \eqref{eq:3.2.1} that
$$g(n)=\sum_{\substack{2i+j=n\\ i,j\ge 0}}p(i)a(j)=\frac{1}{4}\sum_{\substack{2i+j=n\\ i,j\ge 0}}p(i)r(8j+2).$$
Note that $p(n)$ and $r(n)$ are computable by \textit{Mathematica} via functions \texttt{PartitionsP} and \texttt{SquaresR}, respectively. Thus we can complete our verification with much less time. In fact, with the help of \textit{Mathematica}, we see that \eqref{eq:g125.1} and \eqref{eq:g125.2} hold up to the bound $\lfloor v \rfloor=84$, and thus they hold for all $n\ge 0$ by Lemma \ref{le:01}. This completes our proof of Theorem \ref{th:01}. We also end our proof of Theorem \ref{th:02} by a similar verification.

\subsection*{Acknowledgments}

The author thanks S. Radu for interpreting the definition of $\Delta^*$.

\bibliographystyle{amsplain}

\end{document}